\documentclass[12pt,a4paper]{article}
\usepackage[utf8]{inputenc}
\usepackage{fullpage}
\usepackage{float}
\usepackage{amsmath}
\usepackage{amssymb}
\usepackage{stmaryrd}
\usepackage{graphicx}
\usepackage{subcaption}
\usepackage{appendix}
\usepackage{multirow}
\usepackage{physics}
\usepackage{color}
\usepackage{bm}
\usepackage{url}
\usepackage{comment}
\usepackage{hyperref}
\usepackage[linesnumbered,lined,ruled]{algorithm2e}

\usepackage{cite}

\makeatletter
\renewcommand{\@biblabel}[1]{(#1)}
\makeatother

\title{Fast calibration of the LIBOR Market Model with Stochastic Volatility based on analytical gradient}
\author{Herv\'{e} Andres${^1}$, Pierre-Edouard Arrouy$^{1}$, Paul Bonnefoy$^{1}$, \\ Alexandre Boumezoued$^{1}$ and Sophian Mehalla$^{1,2}$ \\ \\
$^1$Milliman R\&D \\
{\small 14 Avenue de la Grande Arm\'ee}\\
{\small Paris, France} \\
$^2$CERMICS \\
{\small 6 - 8 Avenue Blaise Pascal} \\
{\small Cit\'e Descartes – Champs sur Marne} \\
{\small Marne la Vallée, France}
}

\date{\today}

\newcommand{\R}{\mathcal{R}_{m,n}}
\newcommand{\E}{\mathbb{E}}
\newcommand{\pder}[2]{\dfrac{\partial#1}{\partial#2}}
\newcommand{\BTheta}{\mathbf{\Theta}}
\newcommand{\J}{\mathbf{J}}
\newcommand{\f}{\mathbf{f}}
\newcommand{\g}{\mathbf{g}}
\newcommand{\Bd}{\mathbf{d}}
\newcommand{\At}{A(\tau,z)}
\newcommand{\Atj}{A(\tau_j,z)}
\newcommand{\Bt}{B(\tau,z)}
\newcommand{\Btj}{B(\tau_j,z)}

\DeclareMathOperator*{\argmin}{arg\,min}
\newtheorem{prop}{Proposition}
\newtheorem{remark}{Remark}
\begin{document}

\maketitle

\begin{abstract}
We propose to take advantage of the common knowledge of the characteristic function of the swap rate process as modelled in the LIBOR Market Model with  Stochastic Volatility and Displaced Diffusion (DDSVLMM) to derive analytical expressions of the gradient of swaptions prices with respect to the model parameters. We use this result to derive an efficient calibration method for the DDSVLMM using gradient-based optimization algorithms.  

Our study relies on and extends the work by \cite{CUI2017} that developed the analytical gradient for fast calibration of the Heston model, based on an alternative formulation of the Heston moment generating function proposed by \cite{DELBANO2010}.
 
Our main conclusion is that the analytical gradient-based calibration is highly competitive for the DDSVLMM, as it significantly limits the number of steps in the optimization algorithm while improving its accuracy. The efficiency of this novel approach is compared to classical standard optimization procedures. 
\end{abstract}

\section{Introduction}
\label{section_introduction}

To speed up the calibration procedure of the DDSVLMM, two main strategies can be considered. 

First, reduce the computational time required for numerical calculation of one (or multiple) swaption prices. In \cite{WU2006}, pricing under the SVLMM is performed  based on both the classical \cite{HESTON1993} numerical integration method and the famous Fast Fourier Transform (FFT) approach of \cite{CARR1999}, which has become a standard for option valuation for models with known characteristic function, as it is particularly the case for affine diffusion processes. As an alternative to moment generating function calculation, \cite{DEVINEAU2020} have shown the efficiency of using Edgeworth and Gram-Charlier expansions in the calibration of the DDSVLMM. Further work on this type of strategies that reduce the computational cost of the objective function is provided in a companion paper, see \cite{mehalla2020jacobi}.

Note that a typical number of parameters for a DDSVLMM calibration is 8 or 9, depending on whether the displacement parameter is included in the calibration process or not, whereas a standard Heston model calibration involves 5 parameters. In such calibration problem, the regularity of the objective function is unknown. This latter issue can be tackled using optimization algorithms unconcerned about the regularity of the target function; it has been done in \cite{DEVINEAU2020} in which authors used the Nelder-Mead method.

In this paper, we develop a second strategy which consists in decreasing the number of objective function calls required by the optimization algorithm. In this context, the use of gradient-based algorithms is of interest. This can be done by relying on either numerical or analytical gradient computation, the latter being particularly efficient in terms of accuracy and computational cost. A central reference for the present paper is \cite{CUI2017}, that developed the analytical gradient for fast calibration of the Heston model, based on an alternative formulation of the Heston moment generating function proposed by \cite{DELBANO2010}.
Alternatives focusing on `learning' the gradient, without analytical calculation, have been proposed by \cite{LIU2019} based on Artificial Neural Networks. Apart from the Heston model, analytical gradient-based methods have also been developed by \cite{GUDMUNDSSON2019} for the calibration of the so-called 3/2 model. 

By adapting the approach of \cite{CUI2017}, the present paper provides a numerically stable analytical gradient for the DDSVLMM and uses it in gradient-based optimization routines. 
Note that the work of \cite{CUI2017} focused on the Heston model, did not make use of real market data and focused on the Levenberg-Marquardt optimization algorithm. 
We specify the methodology in the context of DDSVLMM swaption modelling and we demonstrate the efficiency of our approach by inputting the analytical gradient into both the Levenberg-Marquardt algorithm (denoted by LM in  the following, see \cite{MARQUARDT1963}) and the Broyden–Fletcher–Goldfarb–Shanno algorithm (denoted by BFGS, see \cite{BYRD1995}) using real market swaption data. 

In addition, we also consider ensuring that the Feller condition is satisfied by using an alternative version of the Levenberg-Marquardt algorithm including constraints; note that although the Feller condition is not targeted in \cite{CUI2017}, in our context this condition is of interest for further uses as it preserves the stability of some numerical discretization schemes. As an example, such Feller condition ensures strong convergence of order 1/2 of discretization schemes as discussed in \cite{ALFONSI2015} (Chapter 3) and references therein.

It is worth mentioning that in our experiment, the computation of the objective function based on the moment generating function representation proposed by \cite{ALBRECHER2007} is in average slightly faster compared to that of \cite{DELBANO2010}. This leads us to consider an `optimal' optimization routine made of 1) the analytical gradient as an input of the gradient-based method (LM or BFGS), based on the differentiation of the moment generating function as provided in \cite{DELBANO2010} and 2) objective function calls during the optimization procedure that rely on the numerical evaluation of the moment generating function representation of \cite{ALBRECHER2007}.

Finally, the efficiency of our approach is compared to the following methods: the classical Heston-type pricing method \cite{HESTON1993} or the Edgeworth expansion of swaption prices as developed in \cite{DEVINEAU2020}, both combined with the Nelder-Mead algorithm as well as the Levenberg-Marquart algorithm with numerical gradient. As a main result, this paper shows that the analytical gradient method is highly competitive both from a computational standpoint and calibration accuracy, as it achieves to significantly limit the number of steps in the optimization algorithm while still offering accurate data replication. 
Our calibration experiments on real data also show that expansion methodologies that reduce the computational cost of an objective function call (see for instance \cite{DEVINEAU2020} and \cite{mehalla2020jacobi}) is a complementary efficient alternative; the combination of the two strategies (expansion approach and computation of the related analytical gradient) appears as a promising direction which is left for further research.

The remainder of this paper is structured as follows. In Section \ref{section_model}, we recall the swap rate (approximated) dynamics under the DDSVLMM, as well as the corresponding moment generating function. Section \ref{section_gradient} presents the alternative representation of the moment generating function we propose along with the analytical gradient calculation of the objective function it implies; the optimization algorithm used is also detailed. Section \ref{section_results} illustrates the efficiency of the proposed calibration method and compares it to the classical alternative methods listed above.

\paragraph*{Notations}
In this work, we consider a probability space $(\Omega, \mathcal{F}, \mathbb{P})$ equipped with a filtration $(\mathcal{F}_t)_{t \geq 0}$ satisfying usual conditions. In our financial context, $\mathbb{P}$ stands for the historical probability measure while $(\mathcal{F}_t)_{t \geq 0}$ represents accumulated market information. Bold font will be employed for denoting either vectors as $\bm{u}$ or matrices using capital letters as $\bm{M}$. The canonical scalar product of two vectors will be denoted by $\bm{u} \cdot \bm{v}$; $\| \bm{u} \|$ will stand for the ($\mathcal{L}^2$-) norm induced by the scalar product: $\| \bm{u} \| = \sqrt{\bm{u} \cdot \bm{u}}$.  

\section{Swaption pricing in the DDSVLMM}
\label{section_model}

The LIBOR Market Model relies on the modelling of the forward rates which are quantities directly observable on the interest rates market. Let $P(t,T)$ be the time-$t$ price of a Zero-Coupon bond maturing at time $T>t$ with par value 1. Let us consider a finite tenor structure $T_0 < T_1 < \dots < T_K <\infty $. For a given $j \in \llbracket 0~,~ K-1 \rrbracket$, the simply compounded forward lending rate, seen at time $t \leq T_j$ and prevailing over the period $[T_j, T_{j+1}[$, is denoted by $F_j(t)$ and an arbitrage-free argument \cite{BRIGO2007} shows that it can be expressed in terms of Zero-Coupon prices as:
\begin{equation}
	F_j(t) = \frac{1}{\Delta T_j} \left( \frac{P(t, T_j)}{P(t, T_{j+1})} - 1\right),
	\label{fwd_dyn}
\end{equation}
with $\Delta T_j := T_{j+1} - T_j$. To account for the modelling of negative forward rates, a displacement coefficient $\delta \geq 0$, often called \textit{shift}, is introduced so that in this framework, the modelled quantities are the shifted-rates $F_j(t) + \delta\text{, } j \in \{0, \dots, K \}$. These are assumed to stay non-negative and a log-normal setting can be prescribed for the modelling. 

Let us consider the spot LIBOR measure $\mathbb{Q}$ - sometimes assimilated to the Risk-Neutral measure - associated with the \text{numéraire} $B(t) = \frac{P(t,T_{m(t)})}{\Pi_{i=0}^{m(t)-1}P(T_i,T_{i+1})}$, where $m(t) = \inf\{0 \leq j \leq K-1 \text{ : } t\le T_j \}$ is the index of the first forward rate that has not expired by $t$. For each $j \in \llbracket 0~,~ K-1 \rrbracket$, we assume the existence of a $N_f$-dimensional deterministic function $t \in \mathbb{R}_+ \longmapsto \bm{\gamma_j}(t)$ such that under $\mathbb{Q}$, the dynamics of the shifted rates writes:
\begin{equation}
\label{lognormal_fwd_rates}
\dd F_j(t) = \sqrt{V_t} (F_j(t) + \delta) \bm{\gamma_j}(t) \cdot \Big( -\bm{ \sigma_{j+1} }(t) \dd t + \dd \bm{W}^*_t \Big),
\end{equation}
where $(\bm{W}^*_t)_{t \geq 0}$ is a $N_f$-dimensional Brownian motion (with independent components) under $\mathbb{Q}$, the function $ \bm{\sigma_{j+1}}(t) := -\sum_{k=m(t)}^j\frac{\Delta T_k(F_k(t)+\delta)}{1+\Delta T_k F_k(t)}\bm{\gamma_j}(t) $ has been introduced and where the process $(V_t)_{t \geq 0}$ is a stochastic process allowing to replicate some market data features (smile, skew). Namely the volatility process lies in the family of Cox-Ingersoll-Ross processes under the Risk-Neutral measure, assuming then the following dynamics: 
\begin{equation}
\label{vol_dyn}
\dd V_t = \kappa ( \theta - V_t) \dd t + \epsilon \sqrt{V_t} \dd W_t
\end{equation}
with $(W_t)_{t \geq 0}$ a Brownian motion and $\kappa, \theta, \epsilon$ three non-negative parameters. The Feller condition $2\kappa\theta \ge \epsilon^2$ ensures that the process remains almost surely non-negative at any time $t$ as long as $V_0 > 0$. 
Finally, the correlation structure between the forward rates and the volatility factor is captured through a time-dependent function $t \in \mathbb{R}_+ \longmapsto \rho_j(t)$ satisfying:
\begin{equation*}
\label{correlation_structure}
\E\left[\left(\frac{ \bm{\gamma_j}(t)}{\| \bm{\gamma_j}(t)\|}\cdot \dd \bm{W}^*_t \right) \dd W_t \right] =: \rho_j(t)\dd t.
\end{equation*}
The specification of the parametric functions $ \bm{\gamma_j}$ and $\rho_j$ are detailed in Section~\ref{calibration_formulation}.

In the remaining of the paper, we will denote by $\BTheta$ the vector of parameters that are to be estimated: it includes the parametric specifications of the maps $\bm{\gamma_j}$ and $\rho_j$ and the volatility parameters $\kappa$, $\theta$ and $\epsilon$. We chose to not include the displaced coefficient $\delta$ in the set of parameters to be estimated and to fix it at an arbitrary value for our experiment.  However, the developed method can be extended to the case when $\delta$ is considered as a parameter to be calibrated. 

\subsection{Swap rate dynamics} 
The time-$t$ value of the (forward) swap rate, $t \leq T_m$, prevailing over the time interval $[T_m, T_n]$ ($m \leq n \leq K$) is given by an arbitrage-free argument again:
\begin{equation*}
\R(t) = \frac{P(t,T_m)-P(t,T_n)}{B^S(t)}
\end{equation*}
where $B^S(t) = \sum_{j=m}^{n-1}\Delta T_j P(t,T_{j+1})$ is the annuity of the considered swap (which stricly depends on $m$ and $n$ although we omit the notation for simplicity). For similar reasons as those previously mentioned, we are lead to consider \textit{shifted} swap rate $\R(t) + \delta$. The \textit{shifted} swap rate can be expressed as a (stochastic) function of the \textit{shifted} forward rates involved between $T_m$ and $T_n$:
\begin{eqnarray*}
	\R(t) + \delta = \sum_{j=m}^{n-1} \alpha_j(t) \big( F_j(t) + \delta \big)
\end{eqnarray*}
with $\alpha_j(t) = \frac{\Delta T_j P(t,T_{j+1})}{B^S(t)}$. Observe that this stochastic weights $\alpha_j$ are themselves functions of the shifted forward rates. Then, the swap rate dynamics under the so-called forward swap rate measure $\mathbb{Q}^S$, associated with the \text{numéraire} $B^S$, can be deduced from the forward rates dynamics (\ref{lognormal_fwd_rates}) using Ito's lemma, see \cite{WU2006}:
\begin{equation*}
\dd \R(t) = \sqrt{V_t}\sum_{j=m}^{n-1}\pder{\R(t)}{F_j(t)} (F_j(t)+\delta) \bm{\gamma_j}(t)\cdot \dd \bm{Z}^S_t, 
\end{equation*}
where $\frac{\partial \R(t)}{\partial F_j(t)} = \alpha_j(t) + \frac{\Delta T_j}{1+\Delta T_j F_j(t)}\sum_{k=m}^{j-1}\alpha_k(t) \big( F_k(t)-\R(t) \big)$ and $(\bm{Z}^S_t)_{t \geq 0}$ is a $N_f$-dimensional Brownian motion. 
 Moreover under $\mathbb{Q}^S$, the dynamics of the volatility factor writes:
\begin{equation*}
\dd V_t = \kappa\left(\theta - \tilde{\xi}^S(t)V_t\right)\dd t +\epsilon\sqrt{V_t}\dd Z_t^S,
\end{equation*}
where $(Z^S_t){t \geq 0}$ is a one-dimensional Brownian motion and \\
$\tilde{\xi}^S(t) := 1 +\frac{\epsilon}{\kappa}\sum_{j=m}^{n-1}\alpha_j(t)\sum_{k=m(t)}^j\frac{\Delta T_k (F_k(t)+\delta)}{1+\Delta T_k F_k(t)}\rho_k(t)\|\bm{\gamma_k}(t)\|$ appeared when applying the Girsanov's theorem. \\

As it stands, the model is too complex to be calibrated. We resort then to the so-called \textit{freezing} technique which relies on the assumption of low variability of some stochastic quantities. We obtain then the following approximated dynamics under $\mathbb{Q}^S$, in which $R_{m,n}(t)$ denotes the approximated shifted swap rate:
\begin{equation}
%\left\{
\label{final_model}
\begin{array}{rcll}
    \dd R_{m,n}(t) & = & \sqrt{V(t)} (R_{m,n}(t)+\delta) \displaystyle\sum_{j=m}^{n-1}\omega_j(0) \bm{\gamma_j}(t)\cdot \dd \bm{Z}^S_t,\\
    \dd V_t & = & \kappa \left(\theta - \tilde{\xi}_0^S(t)V_t\right) \dd t +\epsilon\sqrt{V_t}\dd Z_t^S, 
\end{array}
%\right.
\end{equation}
where $\omega_j(0) := \frac{\partial R_{m,n}(0)}{\partial F_j(0)} \frac{F_j(0)+\delta}{R_{m,n}(0)+\delta} $ and $\tilde{\xi}^S_0(t) := 1 +\frac{\epsilon}{\kappa}\sum_{j=m}^{n-1}\alpha_j(0)\sum_{k=m(t)}^j\frac{\Delta T_k (F_k(0)+\delta)\rho_k(t)\|\bm{\gamma_k}(t)\|}{1+\Delta T_k F_k(0)}$.
Note that (\ref{final_model}) is an Heston-type model (with time-dependent coefficient) that allows to take advantage of well-known pricing method based on the analytical form of the characteristic function of the logarithm of $R_{m,n}$ (the method is detailed in the following section). Indeed, when considering the dynamics of the log-shifted swap rate, one recover an affine diffusion. As a result the moment generating function of $\ln \big( R_{m,n}(t) + \delta \big)$ is explicitly known after solving some Riccati equations.

\subsection{Swaption price}
The spot price of a payer swaption contract with strike $K$ expresses as the following expectation associated to $\mathbb{Q}^S$:
\begin{eqnarray}
\label{price_swp}
PS(\BTheta;0,K) = B^S(0)\E^S[\max(R_{m,n}(T_m)-K,0)],
\end{eqnarray}
where $R_{m,n}(T_m)$ is modelled thanks to (\ref{final_model}). Note that the swaption price is expressed here as a function of the approximated swap rate, so that the swaption price is actually an approximated swaption price. Observe also that we omit the dependence of the swaption price to the maturity and the tenor. It is straightforward to rewrite the swaption price as: 
\begin{equation}
PS(\BTheta;0,K) = B^S(0)\left((R_{m,n}(0)+\delta)P_1(\BTheta;0,K) -(K+\delta)P_2(\BTheta;0,K) \right),
\label{swpt_price}
\end{equation}
where
\begin{align*}
P_1(\BTheta;0,K) &= \E^S\left[e^{\ln{\frac{R_{m,n}(T_m) + \delta}{R_{m,n}(0)+\delta}}}\mathbf{1}_{R_{m,n}(T_m) \ge K} \right], \\
P_2(\BTheta;0,K) &= \E^S\left[\mathbf{1}_{R_{m,n}(T_m) \ge K} \right].
\end{align*}

Let us denote by $(\mathcal{F}_t)_{t\ge 0}$ the filtration generated by the Brownian motions $(\bm{Z^S}_t, Z_t^S)_{t \ge 0}$ and by $\psi$ the moment generating function of the variable $X_{m,n}(T_m) := \ln \bigg( \frac{R_{m,n}(T_m) + \delta}{R_{m,n}(0)+\delta} \bigg) $, defined for $z\in \mathcal{D} \subset \mathbb{C}$ by 
\begin{equation}
\label{MGF}
\psi(\BTheta;X_{m,n}(t),V_t,t;z) = \E^S[e^{zX_{m,n}(T_m)}|\mathcal{F}_t],
\end{equation}
where $\mathcal{D}$ is the domain of convergence of the moment generating function. Let $\varphi$ be the characteristic function of $X_{m,n}(T_m)$ defined for $ z \in \{ z' \in \mathbb{C}: \text{Re}(z') \in \mathcal{D} \}$ by
\begin{equation*}
\label{Characteristic_fun}
\varphi(\BTheta;z) = \psi(\BTheta; X_{m,n}(0), V_0,0;z),
\end{equation*}
where $\text{Re}(\cdot)$ is the real part. In the following, $i$ denotes the imaginary unit satisfying $i^2 = -1$. The two expectations appearing in the swaption price (\ref{swpt_price}) can be computed using the characteristic function $\varphi$ in the following way (we refer the reader to \cite{DUFFIE2000} for a detailed justification): 
\begin{equation*}
\def\arraystretch{3}
\begin{array}{rcl}
P_1(\BTheta;0,K) &=& \dfrac{1}{2}+\dfrac{1}{\pi}\displaystyle\int_0^{+\infty}\Re\left(\dfrac{e^{-iu\ln{\frac{K+\delta}{R_{m,n}(0)+\delta}}}\varphi(\BTheta;u-i)}{iu}\right) \dd u, \\
P_2(\BTheta;0,K) &=& \dfrac{1}{2}+\dfrac{1}{\pi}\displaystyle\int_0^{+\infty}\Re\left(\dfrac{e^{-iu\ln{\frac{K+\delta}{R_{m,n}(0)+\delta}}}\varphi(\BTheta;u)}{iu}\right) \dd u. 
\end{array}
\end{equation*}
As defined in (\ref{MGF}), the moment generating function is a martingale and thus using Ito's lemma and making the drift term zero, we get that $\psi$ is solution of the following Kolmogorov backward equation:
\begin{equation} 
\label{Kolmogorov_eq}
\pder{\psi}{t}+\kappa(\theta- \xi(t) v)\pder{\psi}{v}-\dfrac{1}{2}\lambda^2(t) v \pder{\psi}{x}+\dfrac{1}{2}\epsilon^2 v \frac{\partial^2\psi}{\partial v^2}+\epsilon\tilde{\rho}(t)\lambda(t) v \frac{\partial^2 \psi}{\partial v \partial x}+\dfrac{1}{2}\lambda^2(t) v \frac{\partial^2 \psi}{\partial x^2} = 0,
\end{equation}
with terminal condition 
\begin{equation*}
\psi(\BTheta;x,v,T_m;z) = e^{zx}.
\end{equation*}
The following quantities have been introduced in the equation above:
\begin{equation*}
\xi(t) := \tilde{\xi}^S_0(t),\quad \lambda(t) := \left\|\displaystyle\sum_{j=m}^{n-1}\omega_j(0) \bm{\gamma_j}(t) \right\|,\quad
\tilde{\rho}(t) := \dfrac{1}{\lambda(t)}\displaystyle\sum_{j=m}^{n-1}\omega_j(0)\|\bm{\gamma_j}(t)\|\rho_j(t).
\end{equation*}
Following Heston's approach developed in \cite{HESTON1993} and \cite{WU2006}, one can look for a solution of (\ref{Kolmogorov_eq}) having a separable form $\psi(\BTheta;x,V,t;z) = e^{\At+\Bt V+zx}$ where $\tau = T_m -t$ (the dependence of $A$ and $B$ on $\BTheta$ is omitted for the sake of simplicity). The problem reduces then to the resolution of some Riccati equations (see \cite{WU2006} for the details), which allows to get analytical closed-form expressions of $A$ and $B$ under the common assumption that $\lambda$ and $\tilde{\rho}$ are piecewise constant functions on the grid $(\tau_j,\tau_{j+1}]$ where $\tau_j = T_m - T_{m-j}$, $j = 0, \cdots, m-1$. The recursive computation of $A$ and $B$ is done as follows: for each $j\in \{0,1,\dots,m-1\}$, for each $\tau$ in $(\tau_j , \tau_{j+1}]$,
\begin{equation*}
\left\{
\begin{array}{rcll}
     \At & = & \Atj + \dfrac{\kappa \theta}{\epsilon^2}\left[(\mu+\nu)(\tau-\tau_j)-2\ln\dfrac{1-g_je^{\nu(\tau-\tau_j)}}{1-g_j} \right],& \\
     \Bt & = & \Btj + \dfrac{(\mu+\nu-\epsilon^2\Btj)(1-e^{\nu(\tau-\tau_j)})}{\epsilon^2 (1-g_je^{\nu(\tau-\tau_j)})},& \ \\
\end{array}
\right.
\end{equation*}
with initial condition $A(0,z) = B(0,z) = 0$ and
\begin{equation*}
\mu = \kappa\xi(\tau_j) - \tilde{\rho}(\tau_j)\epsilon\lambda(\tau_j) z,\quad \nu = \sqrt{\mu^2-\lambda^2(\tau_j)\epsilon^2(z^2-z)},\quad g_j=\dfrac{\mu+\nu-\epsilon^2\Btj}{\mu-\nu-\epsilon^2\Btj}.
\end{equation*}
For the sake of simplicity, we omitted the time dependency of $\mu$, $\nu$ and $g_j$.

\section{Analytical gradient calculation and optimization routines}
\label{section_gradient}
\subsection{Analytic characteristic function gradient}
As pointed out by \cite{ALBRECHER2007}, the representation of the characteristic function by \cite{HESTON1993} suffers from numerical discontinuities for long maturities. To overcome this difficulty, they proposed an equivalent representation that is continuous. Another alternative formulation of the Heston moment generating function has been developed by \cite{DELBANO2010} that is also continuous. This alternative formulation is easier to differentiate and thus well suited for gradient calculation, as derived by \cite{CUI2017} in the context of the Heston model.

For our purpose, we rely on the following modifications of the functions $A$ and $B$: for $\tau_j \le \tau < \tau_{j+1},\ j=0,1,\dots,m-1$
\begin{equation}
\label{char_fun_cui}
    \left\{
\begin{array}{rcll}
     A(\tau,z)&=& A(\tau_j,z) -  \dfrac{\kappa\theta\tilde{\rho}\lambda z(\tau-\tau_j)}{\epsilon} + \dfrac{2\kappa \theta}{\epsilon^2}D_j(\tau)&, \\
     B(\tau,z)&=& B(\tau_j,z) -\dfrac{A_j(\tau)}{V_0},\\
\end{array}
\right.
\end{equation}
where:
\begin{align*}
    E_j(\tau) &= \nu+\mu-\epsilon^2B(\tau_j,z)+ (\nu-\mu+\epsilon^2B(\tau_j,z))e^{-\nu(\tau-\tau_j)},  \\
    D_j(\tau) &= \ln{\dfrac{\nu}{V_0}}+ \dfrac{\kappa\xi - \nu}{2}(\tau-\tau_j) - \ln{\dfrac{E_j(\tau)}{2V_0}},  \\
    A^1_j(\tau) &= \left[B(\tau_j,z)(2\mu-\epsilon^2B(\tau_j,z)) + \lambda^2(z-z^2) \right]\tanh{\dfrac{\nu(\tau-\tau_j)}{2}}, \\
    A^2_j(\tau) &= \dfrac{\nu}{V_0}+\dfrac{\mu-\epsilon^2B(\tau_j,z)}{V_0}\tanh{\dfrac{\nu(\tau-\tau_j)}{2}},\\
    A_j(\tau) &= \dfrac{A^1_j(\tau)}{A^2_j(\tau)}.
\end{align*}

Note that these formulas have been adapted from \cite{CUI2017} by using hyperbolic tangent functions in order to avoid numerical instabilities.\\

We now state our key result.
\begin{prop}
The gradient of the swaption price under the approximate Heston dynamics (\ref{final_model}) is given by:

\begin{align*}
 \nabla PS(\BTheta;0,K) = B^S(0)\left((R_{m,n}(0)+\delta)\nabla P_1(\BTheta;0,K) -(K+\delta)\nabla P_2(\BTheta;0,K) \right)
\end{align*}
with
\begin{align*}
	\nabla P_1(\BTheta; 0,K) &= \dfrac{1}{\pi}\int_0^{+\infty}Re\left(\dfrac{e^{-iu\ln{\frac{K+\delta}{R_{m,n}(0)+\delta}}}\nabla\varphi(\BTheta;u-i)}{iu} \right)du, \\
	\nabla P_2(\BTheta; 0,K) &= \dfrac{1}{\pi}\int_0^{+\infty}Re\left(\dfrac{e^{-iu\ln{\frac{K+\delta}{R_{m,n}(0)+\delta}}}\nabla\varphi(\BTheta;u)}{iu} \right)du,
\end{align*}
and 
\begin{align*}
    \nabla\varphi(\BTheta;u) &= \varphi(\BTheta;u)\bm{\chi}(\BTheta;u),
\end{align*}
where $\bm{\chi}(\BTheta;u)$ is the gradient vector (which components are the partial derivatives of the characteristic function with respect to each parameter), detailed in Appendix~\ref{apendix_gradient}.
\label{main_prop}
\end{prop}
\begin{remark}
    The particular form of the analytical gradient allows one to compute $\varphi(\BTheta;u)$ only once. 
\end{remark}
\subsection{DDSVLMM calibration problem formulation}
\label{calibration_formulation}
The calibration problem amounts to find the model parameters that allow to best replicate market data. In this paper, we choose to calibrate on market swaptions prices rather than on implied volatilities since we derived the analytical gradient of the swaption price. \\ 

Let us consider the following standard parametrizations for the maps $\bm{\gamma_j}$ and $\rho_j$ introduced in the section \ref{section_model}: for $t \in [T_k, T_{k+1})$,
\begin{eqnarray*}
  & & \bm{\gamma_j}(t) \equiv \bm{\gamma_j}(T_k) = g(T_j-T_k)\bm{\beta_{j-k+1}}, \\
  & & \rho_j(t) \equiv \rho_j(T_k) = \frac{\rho}{\sqrt{N_f}}\frac{1}{\|\bm{\gamma_j}(T_k)\|} \sum_{p=1}^{N_f}\bm{\gamma_j}^{(p)}(T_k),
\end{eqnarray*}
where $g(u) = (a+bu)e^{-cu}+d$ with $a$, $b$, $c$ and $d$ non-negative parameters to calibrate, $\bm{\beta_{j-k+1}}$ is a $N_f$-dimensional unit vector representing the inter-forward correlation structure and $\rho$ is a correlation parameter to calibrate in $(-1,1)$. Given this parametrization, one can compute the derivatives of $g$ (which is the norm of $\bm{\gamma_j}$ by construction)  with respect to $a$, $b$, $c$ and $d$ as follows, for $u \ge 0$:
\begin{equation*}
    \def\arraystretch{2.2}
\begin{array}{rcl}
    \pder{g(u)}{a}  &=& e^{-cu},\\
    \pder{g(u)}{b}  &=& u\pder{g(u)}{a},\\
    \pder{g(u)}{c}  &=& -(a+bu)\pder{g(u)}{b}, \\
    \pder{g(u)}{d}  &=& 1.\\
\end{array}
\label{d_g_x}
\end{equation*}
Finally, the set of parameters to be calibrated writes $\BTheta = (a,b,c,d,\kappa,\theta,\epsilon,\rho)$. \\

Assume that we have a set $\left(PS^{mkt}_{m,n}(0,K_j)\right)_{j\in J,m \in M, n \in N}$ of market swaptions prices of different strikes, maturities and tenors. For the same set of strikes, maturities and tenors, we denote the model prices computed as in (\ref{swpt_price}) by $\left(PS_{m,n}(0,K_j)\right)_{j\in J,m \in M, n \in N}$. We formulate the calibration problem as an inverse least squares minimization problem with bound constraints and with the Feller condition as an inequality constraint in the following way: \\
\begin{equation}
    \begin{array}{rrcl}
        \displaystyle\argmin_{\BTheta} & \multicolumn{3}{c}{\displaystyle\frac{1}{2W}\displaystyle\sum_{(j,m,n)\in J\times M\times N}w_{j,m,n}\left(\frac{PS_{m,n}(\BTheta;0,K_j)-PS^{mkt}_{m,n}(0,K_j)}{PS^{mkt}_{m,n}(0,K_j)} \right)^2} \\  
         \text{such that } & 2\kappa\theta &\ge& \epsilon^2 \\
            &(a,b,c,d,\kappa,\theta,\epsilon,\rho) &\in& (\mathbb{R}_+)^4\times(\mathbb{R}_+^*)^3 \times (-1;1)\\
    \end{array}
    \label{opti_problem}
\end{equation}
where $w_{j,m,n}$ are fixed positive weights associated to each swaption and \\ $W = \displaystyle\sum_{(j,m,n)\in J\times M\times N} w_{j,m,n}$. In the following, $F$ stands for the objective function in the previous optimization problem and $\f$ stands for the vector of residuals, defined by:
\begin{equation*}
    \f(\BTheta) = \left[\sqrt{\frac{w_{j,m,n}}{W}}\frac{PS_{m,n}(\BTheta;0,K_j)-PS^{mkt}_{m,n}(0,K_j)}{PS^{mkt}_{m,n}(0,K_j)} \right]_{j,m,n}
\end{equation*}
such that we have $F(\BTheta) = \frac{1}{2}\|\f(\BTheta)\|^2$. As mentioned in \cite{GUDMUNDSSON2019}, a regularization term, of the form $\frac{1}{2} \| \Gamma \BTheta\|^2$, can be added to $F$ to promote some solution. As an example, a classical choice is to take $\Gamma = Id$, thus preventing the norm of $\BTheta$ of becoming too large. 

Such an optimization problem can be solved numerically using general optimization methods like the Nelder-Mead algorithm. It is a direct search method, i.e. it does not require any information about the gradient of the objective function. Instead, it relies on the concept of simplex that is a special polytope of $n + 1$ vertices in a $n$-dimensional space. The algorithm starts by evaluating the objective function on each vertex of the simplex. The vertex with the highest value is replaced by a new point where the objective function is lower. The simplex is updated in this manner until the sample standard deviation of the function values on the current simplex falls below some preassigned tolerance. More details on the Nelder-Mead method can be found in \cite{NASH1979}. Note that with this algorithm, one can enforce bound and inequality constraints in (\ref{opti_problem}) by modifying the objective function in the following way:
\begin{align}
\tilde{F}(\BTheta) = 
   \begin{cases}
   F(\BTheta), & \text{if } \BTheta \in (\mathbb{R}_+)^4\times(\mathbb{R}_+^*)^3 \times ]-1;1[ \text{ and } 2\kappa\theta \ge \epsilon^2, \\
   +\infty & \text{otherwise}.
   \end{cases}
\label{F_nm}
\end{align}
Although the Nelder-Mead method has proven to be efficient in some contexts, we observed that it turns out to be very time-consuming in our framework. As a matter of fact, this method requires a lot of evaluations of the objective function and the computation of the swaption prices is very expensive due to the computation of integrals in the complex field and the recursive definition of the characteristic function. This has already been pointed out by \cite{DEVINEAU2020}. Consequently, an optimization algorithm that does not require a lot of objective functions evaluations is preferred in order to achieve a fast calibration. Since we have been able to derive an analytical formula for the gradient of model swaption price, we study gradient-based optimization algorithms. We present two of such methods in the next section. 

\subsection{Calibration using gradient-based algorithms}
In this section, we quickly remind the main ideas behind gradient-based algorithms before presenting more specifically two of these algorithms, namely the BFGS and the LM algorithms. \\

Gradient-based algorithms are iterative methods that start from a given point and proceed by successive adjustments. Each improvement is obtained by moving from the current point along a conveniently chosen descent direction in such a manner that the value of the objective function decreases. We describe in Algorithm \ref{alg:grad_based_algos} the general algorithmic scheme of gradient-based algorithms. More details on gradient-based algorithms can be found in \cite{NOCEDAL2006}.\\

\begin{algorithm}[H] 
\KwIn{Initial guess $\BTheta_0$, objective function $F$, objective function gradient $\nabla F$, tolerance $\epsilon_{tol}$ and maximum number of iterations $k_{max}$}
\Begin{ 
$k \gets 0$ \\
\While{$k \le k_{max}$}{
    Compute $\nabla F(\BTheta_k)$ \\
    \If{$\nabla F(\BTheta_k) \le \epsilon_{tol}$}{
        \textbf{break}
    }
    Compute a descent direction $\bm{d}_k$, generally by using $\nabla F(\BTheta_k)$ \\
    Find $\alpha_k$ such that $F(\BTheta_k + \alpha_k \bm{d}_k)$ is reasonably lower than $F(\BTheta_k)$ (line search) \\
    $\BTheta_{k+1} \gets \BTheta_k + \alpha_k \bm{d}_k$
}
}
\KwOut{Last computed $\BTheta_k$}
\caption{Gradient-based algorithms in pseudo code}
\label{alg:grad_based_algos}
\end{algorithm}

We now take a closer look at the BFGS and LM algorithms. The BFGS algorithm is a quasi-Newton method, which means that it uses an approximation of the inverse of the Hessian matrix instead of the exact inverse used in Newton's method. It is designed for solving all kinds of unconstrained non-linear optimization problems. The descent direction $\bm{d}_k$ is given by
$$
\bm{d}_k = -\mathbf{H}_k\nabla F(\BTheta_k),
$$
where $\mathbf{H}_k$ is an approximation of the inverse of the Hessian matrix that is defined recursively by,
$$
\mathbf{H}_{k+1} = \mathbf{H}_k - \frac{\mathbf{H}_k \mathbf{y}_k \mathbf{y}_k^T\mathbf{H}_k}{\mathbf{y}_k^T \mathbf{H}_k \mathbf{y}_k}+\frac{\mathbf{s}_k \mathbf{s}_k^T}{\mathbf{y}_k^T \mathbf{s}_k}
$$
with initial value $\mathbf{H}_0 = \mathbf{I}$ ($\mathbf{I}$ is the identity matrix), $\mathbf{s}_k = \BTheta_{k+1} - \BTheta_k$ and $\mathbf{y}_k = \nabla F(\BTheta_{k+1}) - \nabla F(\BTheta_k)$. The gradient of the objective function $F$ writes as a function of the gradient of the swaption price $\nabla PS_{m,n}$:
\begin{equation*}
\nabla F(\BTheta) = \frac{1}{W}\sum_{(i,m,n)\in I\times M\times N}w_{j,m,n}\nabla PS_{m,n}(\BTheta;0,K_j)\frac{PS_{m,n}(\BTheta;0,K_j) - PS^{mkt}_{m,n}(0,K_j)}{PS^{mkt}_{m,n}(0,K_j)^2}.
\end{equation*}
The calibration problem (\ref{opti_problem}) is a constrained optimization problem and thus, the BFGS algorithm can not be used. We therefore rely on an extension of the classical BFGS algorithm, known as L-BFGS-B, that has been developed to handle bound constraints \cite{BYRD1995}. However, inequality constraints can not be easily enforced. Consequently, we relax the Feller condition when using this method. 

Furthermore, the Levenberg-Marquardt algorithm is specifically designed for solving non-linear least squares problems, which is exactly the type of problem we are coping with. This algorithm has the particularity of behaving like the steepest descent method when the current point is far (in some sense) of a (the) solution and of behaving like the Gauss-Newton method when the current point is near of a (the) solution. This is achieved by introducing a damping parameter in the expression of the descent direction, as follows
\begin{equation*}
    \Bd_k = -(\J_k^T \J_k +\mu_k \mathbf{I})^{-1}\g_k
\end{equation*}
where $\J_k$ is the Jacobian matrix of $\f$ in $\BTheta_k$, $\g_k=\J_k^T\f(\BTheta_k)$ is the gradient of $F$ in $\BTheta_k$ and $\mu_k$ is the damping parameter. For large values of $\mu_k$ (compared to coefficients of $\J_k^T \J_k$), we have $\Bd_k \simeq -\frac{1}{\mu_k}\g_k $ which corresponds to the descent direction in the steepest descent method. For small values of $\mu_k$, we have $\Bd_k \simeq -(\J_k^T \J_k)^{-1}\g_k$ which corresponds to the descent direction in the Gauss-Newton method. The updating strategy of $\mu_k$ is described in Algorithm \ref{alg:lm_algo} in Appendix~\ref{classic_lm}. There are a plenty of different updating strategies leading to several versions of the Levenberg-Marquardt algorithm. More details on this algorithm can be found in \cite{MADSEN1999}.

\paragraph*{Adding bound constraints to the LM algorithm\\}
As for the BFGS algorithm, the classic LM algorithm does not handle bound and inequality constraints. However, it can be extended to do so.
Bound constraints can be ensured by using a projection of $\BTheta_k$ onto the feasible set. We detail the modifications in Algorithm \ref{alg:lmbc_algo} in Appendix~\ref{extended_lm}. Handling linear inequality constraints requires many modifications so we will not detail them here but the interested reader can find a discussion on this topic in \cite{NOCEDAL2006} (Chapter 15). Observe that the Feller condition is not linear in the parameters of the DDSVLMM. However it can be easily linearized using the following change of variables: 
\begin{equation*}
    \tilde{\kappa} = \ln{\kappa},\quad \tilde{\theta} = \ln{\theta},\quad \tilde{\epsilon} = \ln{\epsilon}.
\end{equation*}
The Feller condition writes in term of the new volatility parameters as:
\begin{equation*}
    \tilde{\kappa} + \tilde{\theta} + \ln{2} \ge 2\tilde{\epsilon}.
\end{equation*}
To account for this change of variables, the gradient of the swaption price has to be modified by replacing the partial derivatives with respect to $\kappa$, $\theta$ and $\epsilon$ by the corresponding derivatives with respect to $\tilde{\kappa}$, $\tilde{\theta}$ and $\tilde{\epsilon}$:
\begin{equation*}
\def\arraystretch{2.2}
\begin{array}{rcl}
\pder{PS}{\tilde{\kappa}}(\BTheta;0,K_j) &=& \kappa\pder{PS}{\kappa}(\BTheta;0,K_j), \\
\pder{PS}{\tilde{\theta}}(\BTheta;0,K_j) &=& \theta\pder{PS}{\theta}(\BTheta;0,K_j), \\
\pder{PS}{\tilde{\epsilon}}(\BTheta;0,K_j) &=& \epsilon\pder{PS}{\epsilon}(\BTheta;0,K_j). \\
\end{array}
\end{equation*}

\section{Calibration results}\label{section_results}
In this section, we present our experimental results for the calibration of the DDSVLMM (\ref{final_model}) using the BFGS and LM algorithms. We first detail the market data to be replicated and discuss some implementation aspects. Then, we compare the BFGS and LM routines with existing calibration methods with regards to the objective function value and to the computational time. \\

For the calibration, we used a set of 280 market EURO and USD swaptions volatilities. For the purpose of the calibration, these volatilities are converted into prices using the Bachelier formula based on a rate curve as used under the Solvency II regulation\footnote{Available at \url{www.eiopa.europa.eu}}. The ATM swaptions maturities and tenors considered range into $\{1,\dots,10,15,20,25,30\}$. For away-from-the-money swaptions, we consider the same range for maturities and focus on a 10-year reference tenor; the strikes (in bps) range into $+/- \{25,50,100\}$. As mentioned previously, the shift $\delta$ is objectified otherwise: we take $\delta = 0.1$. The inter-forward correlation structure, captured by the $\bm{\beta_k}$ parameters, is assessed by an PCA technique we do not detail here. The number $N_f$ of risk factors is set to 2.  \\

We implemented the pricing and gradient functions in C++. We used the R base function \texttt{optim} for the Nelder-Mead and BFGS algorithms and we used the C++ \texttt{LEVMAR} package \cite{LOURAKIS2005} for the Levenberg-Marquardt algorithm. This choice is particularly motivated by the fact that the \texttt{LEVMAR} package implements the extended version of the Levenberg-Marquardt algorithm allowing to handle bound constraints and the extended version allowing to handle both bound and linear inequality constraints. As for the computation of the numerical integral required in the pricing and the gradient functions, we resorted to the Gauss-Laguerre quadrature with 90 nodes. 

\subsection{Methods accuracy}
We compare the BFGS and Levenberg-Marquardt algorithms with existing calibration methods based on the criteria of the objective function value. First, let us introduce the three reference calibration methods used for the purpose of comparison.  \\

The first one is the classical Heston method in which the price is computed through the formula (\ref{swpt_price}) and the optimization is performed via the Nelder-Mead algorithm. We set the maximum number of iterations to 500 and we repeat the optimization 3 times in order to achieve a better convergence.

The second calibration method relies on Edgeworth approximations: it consists in using an approximate swaption price obtained using an Edgeworth expansion of the unknown density of the swap rate $\R$ (see \cite{DEVINEAU2020} for a thorough description of the method). The associated optimization method is the Nelder-Mead algorithm. We use the same parametrization for the Nelder-Mead method as for the classical Heston method. 

The last method is based on the LM algorithm but associated with a numerical gradient estimation instead of using the derived analytical gradient. The price is still computed with pricing formula (\ref{swpt_price}). We use the central difference scheme in order to approximate the gradient as:
\begin{equation*}
    \nabla PS(\BTheta;0,K_j) \simeq \frac{PS(\BTheta + \bm{h};0,K_j,T_m,T_n)-PS(\BTheta-\bm{h};0,K_j)}{2\bm{h}}
\end{equation*}
where $\bm{h} := h \bm{e}$ with $\bm{e}$ a vector whose components are $1$ and $h$ a small quantity. We take $h = 10^{-8}$ and a maximal number of 15 iterations.\\

Let us also present the different parametrizations studied for the BFGS and Levenberg-Marquardt algorithms relying on the analytical gradient as derived in this paper. For the BFGS algorithm, we test one configuration in which the maximum number of iterations is set to 30. For the Levenberg-Marquardt algorithm, we test two configurations having the same tolerance levels $\epsilon_1$, $\epsilon_2$ and $\epsilon_3$ that are set to $10^{-10}$ (see Appendix~\ref{classic_lm}). The two configurations, respectively denoted by LM-BC-15 and LM-BC-30, use the version of the LM algorithm allowing to handle bound constraints only with a maximum number of iterations of 15 and 30 respectively.\\

We summarize the studied methods and their main characteristics in Table \ref{tab:summary_methods}.

\begin{table}[H]
\resizebox{\textwidth}{!}{%
\begin{tabular}{|ccccc|}
\hline
Method name & \begin{tabular}[c]{@{}c@{}}Optimization \\ method\end{tabular} & \begin{tabular}[c]{@{}c@{}}Maximum number\\ of iterations\end{tabular} & \begin{tabular}[c]{@{}c@{}}Ensures Feller\\ condition\end{tabular} & Features \\ \hline
Heston & Nelder-Mead & 500 & Yes & Repeated 3 times \\
Edgeworth & Nelder-Mead & 500 & Yes & \begin{tabular}[c]{@{}c@{}}Repeated 3 times. Uses a different\\ swaption pricing formula.\end{tabular} \\
LM-NUM & Levenberg-Marquardt & 15 & No & Uses the numerical gradient \\
BFGS & L-BFGS-B & 30 & No & Uses the analytical gradient \\
LM-BC-30 & Levenberg-Marquardt & 30 & No & Uses the analytical gradient \\
LM-BC-15 & Levenberg-Marquardt & 15 & No & Uses the analytical gradient \\ \hline
\end{tabular}%
}
\caption{Summary of studied methods}
\label{tab:summary_methods}
\end{table}

Before going further, we precise the bounds, particularly the lower bounds, for the 8 parameters of the DDSVLMM to calibrate. Indeed, we experienced numerical instabilities when some parameters equal zero or are very close to zero. For instance if the speed reversion parameter $\kappa$ or the volatility of volatility $\epsilon$ become almost zero, the behavior of the model is pathologic. Therefore, we give the following lower ($\bm{LB}$) and upper ($\bm{UB}$) bounds for $\BTheta = (a,b,c,d,\kappa,\theta,\epsilon,\rho)$:
\begin{equation*}
\begin{array}{rcl}
    \bm{LB} &:=& (10^{-5},10^{-5},10^{-5},10^{-5},10^{-5},10^{-5},10^{-5},-0.999), \\
     \bm{UB} &:=& (+\infty,+\infty,+\infty,+\infty,+\infty,+\infty,+\infty,0.999).
\end{array}
\end{equation*}

\begin{figure}[H]
    \centering
    \includegraphics[width=\textwidth]{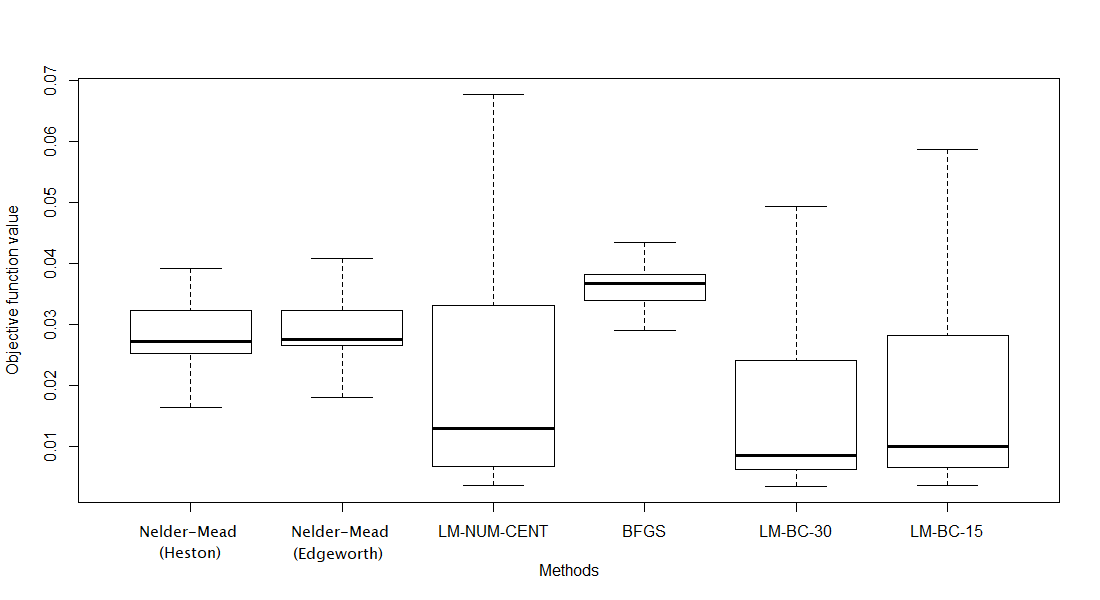}
    \caption{Boxplots of our selected benchmark of calibration methods.}
    \label{fig:boxplot_eur}
\end{figure}

The procedure that has been led in order to compare the various calibration methods is the following: for each set of data, we draw randomly 100 initial parameter starting values between $\bm{LB}$ and $\bm{UB}$ that satisfy the Feller condition and we perform the calibration for each described method starting from each of these initial guess. From this procedure, we retrieve the boxplots of Figure \ref{fig:boxplot_eur} using EURO data, which provide statistics on the objective function value over the 100 calibrations.

Note first that the variance of replication errors when using Nelder-Mead and BFGS algorithms are comparable and significantly lower to that obtained when using Levenberg-Marquardt routines. To this extent, Nelder-Mead and BFGS seem less dependent to initial guess. However, median errors reached by the Levenberg-Marquardt are lower: this can be explained by the fact that this algorithm is particularly suited for least-squares problems. 
	
When working with USD data, for which results are presented in Appendix~\ref{US_data_results}, these conclusions still hold. Nevertheless, Nelder-Mead based techniques as well as BFGS using analytical gradient perform significantly better on USD data. It can be explained by the fact that the Nelder-Mead algorithm is uniquely based on the topology of the graph of the objective function and thus substantially depends on market data in our context. 

Concerning the Levenberg-Marquardt approach more specifically, we first note that increasing the number of iterations  yields a better calibration (in average and in variance), which is an expected behaviour. The improvement is significant using USD data. Using the Levenberg-Marquardt technique coupled with a numerical approximation of the gradient of the target function leads to a rather wide range of objective function values, which illustrates the benefit of using an analytical gradient in comparison. In addition, the median value for this approach is higher to that obtained when using Levenberg-Marquardt optimization with analytical gradient. Therefore, the information conveyed by the analytical Hessian matrix turns out to be valuable in order to stabilize the calibration process and reduce the dependency to the starting point. This conclusion holds for both EURO and USD data.

So far we did not impose the Feller's condition to be satisfied by the outputs of the calibration procedures. The percentages of obtained parameters (over 100 calibrations) that do not satisfy it are given in Table~\ref{tab:feller_condition}. The number of unsatisfied Feller conditions is rather significant especially for the LM algorithm. Note that such condition is always ensured for Nelder-Mead based calibrations as the Feller condition has been imposed as depicted in (\ref{F_nm}). When working with USD data, behaviours of the different algorithms are similar.

\begin{table}[H]
\centering
\begin{tabular}{|cc|}
\hline
Method &  \begin{tabular}[c]{@{}c@{}}Percentage of unsatisfied \\ Feller condition \end{tabular} \\ \hline
Heston & 0 \% \\
Edgeworth & 0 \% \\
LM-NUM & 34 \% \\
BFGS & 13 \%  \\
LM-BC-30 & 33 \% \\ 
LM-BC-15 & 32 \%  \\ \hline
\end{tabular}
\caption{Percentages of unsatisfied Feller condition over the 100 calibrations for each methods}
\label{tab:feller_condition}
\end{table}

\begin{figure}[H]
  \centering
  \includegraphics{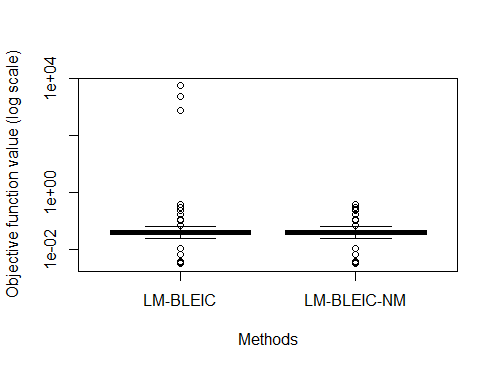}
  \caption{Boxplots for LM-BLEIC and LM-BLEIC-NM (log-scale).}
  \label{fig:boxplot_lmbleic}
\end{figure}

\paragraph*{Numerical results with bound constraints} 
In view of the previous results, we propose to study two other configurations of the Levenberg-Marquardt algorithm. The first one, denoted by LM-BLEIC, uses the version of the Levenberg-Marquardt algorithm allowing to handle both bound constraints and linear inequality constraints. The maximum number of iterations is set to 50. The second one, denoted by LM-BLEIC-NM, uses the same algorithm as the LM-BLEIC configuration but here 200 iterations of the Nelder-Mead method are peformed at the end of the Levenberg-Marquardt algorithm when the latter converged towards a point whose objective function value is greater than a given threshold set to 0.3. We present the boxplots for the LM-BLEIC and LM-BLEIC-NM methods in Figure \ref{fig:boxplot_lmbleic} using EURO data.
We observe that for some initial guesses, the LM-BLEIC method converges towards points at which the objective function takes extremely high values. This is actually due to the fact that the implementation of the Levenberg-Marquardt algorithm handling linear inequality constraints can not ensure that the points stay in the feasible set. As a consequence, some parameters may take negative values leading to numeric precision errors forcing the algorithm to stop prematurely. This is the reason why we introduced the LM-BLEIC-NM. When the objective function value at the exit of the LM-BLEIC method exceeds a given threshold, we perform some Nelder-Mead iterations in order to achieve a better optimum. This method gives us quite satisfying results compared to LM-BLEIC, as far as it allows to get rid of the extreme points. Using the Nelder-Mead allows thus to reduce the variance of the outputs of the calibration procedures but does not significantly reduce the median value of the replication errors which remain close. Similar conclusions hold when comparing those two methods on USD data as depicted in Figure \ref{fig:boxplots_lmbleic_usd} in the Appendix \ref{US_data_results}. 

\subsection{Time efficiency}
Following the same procedure as the one presented in the previous section, we compute the average CPU time, the average number of calls to the objective function and the average duration of a call. For gradient-based methods, we also compute the average number of calls to the gradient and the average duration of a call. The computations were performed on a computer with a 2.6 gigahertz Intel Core i7 processor and 8 gigabytes of RAM. The results are presented in Table \ref{tab:cpu_time_eur} for EURO data and in Appendix \ref{US_time} in Figure \ref{tab:cpu_time_usd} for USD data. Note that when using the Levenberg-Marquardt algorithm, we compute the average call time to the residual function $f$ and its gradient $\nabla f$ instead of the objective function $F$ and its gradient $\nabla F$ because the \texttt{LEVMAR} implementation takes as inputs the functions $f$ and $\nabla f$, see Section \ref{calibration_formulation}. This explains why there is a difference between BFGS and Levenberg-Marquardt methods in terms of average call time to the gradient. Note also that calling the objective function in the Heston method appears to take in average less time than in gradient-based methods since $F$ is actually replaced by $\tilde{F}$ defined in (\ref{F_nm}) for Nelder-Mead methods, which simply returns a large value in all cases where the Feller condition is not satisfied. \\

\begin{table}[H]
\resizebox{\textwidth}{!}{
\begin{tabular}{|cccccc|}
\hline
Method & Average CPU time & \begin{tabular}[c]{@{}c@{}}Average number of \\ calls to $F$/$f$\end{tabular} & \begin{tabular}[c]{@{}c@{}}Average call time to \\ $F$/$f$\end{tabular} & \begin{tabular}[c]{@{}c@{}}Average number of \\ calls to $\nabla F$/$\nabla f$\end{tabular} & \begin{tabular}[c]{@{}c@{}}Average call time to\\ $\nabla F$/$\nabla f$\end{tabular} \\ \hline
Heston & 159.45 s & 1489.20 & 0.11 s & 0 & 0 s \\
Edgeworth & 8.47 s & 1469.46 & 5.55 $\times 10^{-3}$ s & 0 & 0 s \\
LM-NUM-CENT & 40.06 s & 277.36 & 0.14 s & 0 & 0 s \\
BFGS & 36.40 s & 39.34 & 0.14 s & 39.34 & 0.78 s \\
LM-BC-30 & 33.89 s & 104.80 & 0.14 s & 30.00 & 0.63 s \\
LM-BC-15 & 14.94 s & 38.93 & 0.14 s & 15.00 & 0.63 s \\
LM-BLEIC & 44.29 s & 95.76 & 0.14 s & 48.34 & 0.63 s \\
LM-BLEIC-NM & 45.51 s & 102.29 & 0.14 s & 48.34 & 0.63 s \\ \hline
\end{tabular}%
}
\caption{Computational times}
\label{tab:cpu_time_eur}
\end{table}

The gradient-based algorithms (including those using numerical gradient) appear to be much faster than the classical Heston calibration method using the Nelder-Mead algorithm, since they provide reductions in computational time ranging from 71 \% (LM-BLEIC-NM) to 91 \% (LM-BC-15). This gain in time results directly from the massive reduction of the number of calls to the objective function and thus, to the characteristic function. However, these methods are still not as fast as the Edgeworth method which uses a large number of function calls but for which each call is very fast, which was the purpose of the method at its origin, see \cite{DEVINEAU2020}. However one needs to keep in mind that the reduction in computational time achieved by the Edgeworth expansion comes at the cost of a lower fitting accuracy to market data, as pointed out in \cite{DEVINEAU2020}. % and \cite{mehalla2020jacobi}.

As a main conclusion, the use of the analytical gradient rather than the numerical gradient allows to reduce the calibration duration by a factor of 2.7 when looking at LM-BC-15.

Finally, let us observe that the LM algorithm using analytical gradient and including bound constraints only (LM-BC) is faster than the BFGS algorithm. This can be explained by the fact that the Levenberg-Marquardt algorithm takes advantage from the particular shape of the calibration problem, i.e. a least squares optimization problem. \\

To conclude this section, we justify numerically why we use the characteristic function by \cite{ALBRECHER2007} instead of that of \cite{DELBANO2010} in the computation of the swaption price. We compared the average call time of both characteristic function representations over 1000 randomly chosen parameters and observed an average call time to Albrecher's representation 5 \% lower than the average call time to Cui's representation. The time difference between these two representations of the characteristic function is essentially due to the fact that the coefficients $A$ and $B$ in Albrecher's representation can be easily written as function of a few quantities (see e.g. Appendix 5.3 of \cite{DEVINEAU2020}) which allows to perform less computations than for Cui's representation. 

\appendix
%\appendixpage
\section{Gradient expression}
\label{apendix_gradient}
The vector $\bm{\chi}$ in Proposition \ref{main_prop} writes:
\begin{equation*}
    \bm{\chi}(\BTheta;z) := \left[\chi_a(z), \chi_b(z), \chi_c(z), \chi_d(z), \chi_{\kappa}(z), \chi_{\theta}(z), \chi_{\epsilon}(z), \chi_{\rho}(z) \right]^T
\end{equation*}
where $\chi_x$ denotes the partial derivative of $\varphi$ with respect to the parameter $x$. Using that $\chi_x(z) = \frac{\partial \psi(\BTheta; X_{m,n}(0), V_0,0;z)}{\partial x}$ and since $\psi(\BTheta; X_{m,n}(0), V_0,0;z)$ is defined recursively with terminal value $\psi(\BTheta;X_{m,n}(0),V_0,T_m;z)$, one needs to compute $\frac{\partial \psi(\BTheta; X_{m,n}(0), V_0,t;z)}{\partial x}$ on each interval $(\tau_j,\tau_{j+1}]$. We will rather give the partial derivatives of $\psi$ at any time. 

\subsection{Partial derivative of \texorpdfstring{$\psi$}{psi} with respect to \texorpdfstring{$\theta$}{theta}}
Since $A_j(\tau)$ and consequently also $\Bt$ are independant from $\theta$ for all $t$, the partial derivative of $\psi$ with respect to $\theta$ writes:  
$$
    \dfrac{\partial \psi(X_{m,n}(t),V(t),t;z)}{\partial \theta} = \dfrac{\partial \At}{\partial \theta}\psi(X_{m,n}(t),V(t),t;z).
$$
The partial derivative of $\At$ is given by:
$$
\dfrac{\partial \At}{\partial \theta} = \dfrac{\partial \Atj}{\partial \theta} - \dfrac{\kappa \tilde{\rho}\lambda z(\tau-\tau_j)}{\epsilon} + \dfrac{2\kappa}{\epsilon^2}D_j(\tau)
$$
since $D_j(\tau)$ is independant from $\theta$.

\subsection{Partial derivative of \texorpdfstring{$\psi$}{psi} with respect to \texorpdfstring{$\kappa$}{kappa}}
The partial derivative of $\psi$ with respect to $\kappa$ writes: 
$$
    \dfrac{\partial \psi(X_{m,n}(t),V(t),t;z)}{\partial \kappa} = \left[\dfrac{\partial \At}{\partial \kappa} + V_0\dfrac{\partial \Bt}{\partial \kappa} \right]\psi(X_{m,n}(t),V(t),t;z)
$$
with 
$$
\def\arraystretch{2.2}
\begin{array}{rcl}
  \dfrac{\partial \At}{\partial \kappa} &=& \dfrac{\partial \Atj}{\partial \kappa}-\dfrac{\theta\tilde{\rho}\lambda z(\tau-\tau_j)}{\epsilon} + \dfrac{2\theta}{\epsilon^2}D_j(\tau) + \dfrac{2\kappa\theta}{\epsilon^2}\dfrac{\partial D_j(\tau)}{\partial \kappa}, \\
  \dfrac{\partial \Bt}{\partial \kappa} &=& \dfrac{\partial \Btj}{\partial\kappa} - \dfrac{1}{V_0}\dfrac{\partial A_j(\tau)}{\partial \kappa},
\end{array}
$$
where
$$
\def\arraystretch{2.7}
\begin{array}{rcl}
    \dfrac{\partial D_j(\tau)}{\partial \kappa} &=& \dfrac{\mu}{\nu^2} + \dfrac{1}{2}\left(1-\dfrac{\mu}{\nu}\right)(\tau-\tau_j) - \dfrac{1}{E_j(\tau)}\pder{E_j(\tau)}{\kappa}, \\
    \dfrac{1}{E_j(\tau)}\pder{E_j(\tau)}{\kappa} &=& \dfrac{1}{\nu V_0 A_j^2(\tau)}\left[ \mu + \nu\left(1-\epsilon^2\dfrac{\partial \Btj}{\partial \kappa}\right)\tanh{\dfrac{\nu(\tau-\tau_j)}{2}}\right],, \\
    && - \dfrac{1}{E_j(\tau)}(\nu-\mu+\epsilon^2\Btj)\dfrac{\mu}{\nu}(\tau-\tau_j)e^{-\nu(\tau-\tau_j)},\\
    \pder{A_j(\tau)}{\kappa} &=&  \dfrac{1}{\tilde{A}^2_j(\tau)}\pder{\tilde{A}^1_j(\tau)}{\kappa}-\dfrac{A_j(\tau)}{\tilde{A}^2_j(\tau)}\pder{\tilde{A}^2_j(\tau)}{\kappa},\\
      \dfrac{1}{\tilde{A}^2_j(\tau)}\pder{\tilde{A}^1_j(\tau)}{\kappa} &=& \dfrac{1}{A^2_j(\tau)}\bigg[2\tanh{\dfrac{\nu(\tau-\tau_j)}{2}}\left(\mu\pder{\Btj}{\kappa}-\epsilon^2\Btj)\pder{\Btj}{\kappa}+\Btj\right) \\
    && +\dfrac{\mu(\tau-\tau_j)}{2\nu}\left(\Btj(2\mu-\epsilon^2\Btj)+\lambda^2(z-z^2)\right)\bigg],\\
    \dfrac{1}{\tilde{A}^2_j(\tau)}\pder{\tilde{A}^2_j(\tau)}{\kappa}&=& \dfrac{1}{V_0 A^2_j(\tau)}\bigg[\left(1+\dfrac{\mu(\tau-\tau_j)}{2}-\epsilon^2\pder{\Btj}{\kappa}\right)\tanh{\dfrac{\nu(\tau-\tau_j)}{2}}\\
    &&+ \dfrac{\mu}{\nu}\left(\dfrac{\tau-\tau_j}{2}(\mu-\epsilon^2\Btj)+1 \right) \bigg].
\end{array}
$$
Note that in order to make the calculation of the partial derivative of $A_j$ easier, we write $A_j$ as the ratio of $\tilde{A}_j^1$ and 
$\tilde{A}_j^2$ instead of $A_j^1$ and $A_j^2$, where $\tilde{A}_j^1$ and 
$\tilde{A}_j^2$ are defined as: 
$$
\begin{aligned}
\tilde{A}_j^1(\tau) &=& A_j^1(\tau)\cosh{\dfrac{\nu(\tau-\tau_j)}{2}}, \\
\tilde{A}_j^2(\tau) &=& A_j^2(\tau)\cosh{\dfrac{\nu(\tau-\tau_j)}{2}}.
\end{aligned}
$$
This trick will be re-employed for other derivatives. 
\subsection{Partial derivative of \texorpdfstring{$\psi$}{psi} with respect to \texorpdfstring{$\epsilon$}{epsilon}}
The partial derivative of $\psi$ with respect to $\epsilon$ writes: 
$$
	\dfrac{\partial \psi(X_{m,n}(t),V(t),t;z)}{\partial \epsilon} = \left[\dfrac{\partial \At}{\partial \epsilon} + V_0\dfrac{\partial \Bt}{\partial \epsilon} \right]\psi(X_{m,n}(t),V(t),t;z)
$$
with
$$
\def\arraystretch{2.7}
\begin{array}{rcl}
    \pder{\At}{\epsilon} &=& \pder{\Atj}{\epsilon}+ \dfrac{\kappa\theta\tilde{\rho}\lambda z (\tau-\tau_j)}{\epsilon^2} - \dfrac{4\kappa\theta}{\epsilon^3}D_j(\tau) + \dfrac{2\kappa\theta}{\epsilon^2}\pder{D_j(\tau)}{\epsilon},\\
    \pder{\Bt}{\epsilon} &=& \pder{\Btj}{\epsilon} - \dfrac{1}{V_0}\pder{A_j(\tau)}{\epsilon},
\end{array}
$$
where
$$
\def\arraystretch{2.7}
\begin{array}{rcl}
     \pder{\xi}{\epsilon} &=& \dfrac{\xi - 1}{\epsilon}, \\
     \pder{\mu}{\epsilon} &=& \dfrac{\mu-\kappa}{\epsilon}, \\
     \pder{\nu}{\epsilon} &=& \dfrac{\nu^2-\kappa\mu}{\epsilon\nu}, \\
     \pder{D_j(\tau)}{\epsilon} &=& \dfrac{1}{\nu}\pder{\nu}{\epsilon} + \dfrac{1}{2}\left(\kappa\pder{\xi}{\epsilon}-\pder{\nu}{\epsilon}\right)(\tau-\tau_j) - \dfrac{1}{E_j(\tau)}\pder{E_j(\tau)}{\epsilon},\\
     \dfrac{1}{E_j(\tau)}\pder{E_j(\tau)}{\epsilon} &=& \dfrac{1}{V_0A_j^2(\tau)}\left[\pder{\nu}{\epsilon}+\left(\pder{\mu}{\epsilon}-2\epsilon\Btj-\epsilon^2\pder{\Btj}{\epsilon} \right)\tanh{\dfrac{\nu(\tau-\tau_j)}{2}}\right]\\
     &&-\dfrac{1}{E_j(\tau)}(\nu-\mu+\epsilon^2\Btj)\pder{\nu}{\epsilon}(\tau-\tau_j)e^{-\nu(\tau-\tau_j)},\\
     \pder{A_j(\tau)}{\epsilon} &=&  \dfrac{1}{\tilde{A}^2_j(\tau)}\pder{\tilde{A}^1_j(\tau)}{\epsilon}-\dfrac{A_j(\tau)}{\tilde{A}^2_j(\tau)}\pder{\tilde{A}^2_j(\tau)}{\epsilon},\\
         \dfrac{1}{\tilde{A}_j^2(\tau)}\pder{\tilde{A}_j^1(\tau)}{\epsilon} &=& \dfrac{1}{A_j^2(\tau)}\bigg[2\tanh{\dfrac{\nu(\tau-\tau_j)}{2}}\bigg(\mu\pder{\Btj}{\epsilon}-\epsilon^2\Btj\pder{\Btj}{\epsilon}+\pder{\mu}{\epsilon}\Btj \\
         &&-\epsilon(\Btj)^2 \bigg)+\pder{\nu}{\epsilon}\dfrac{\tau-\tau_j}{2}\left(\Btj(2\mu-\epsilon^2\Btj)+\lambda^2(z-z^2) \right) \bigg],\\
     \dfrac{1}{\tilde{A}_j^2}\pder{\tilde{A}_j^2}{\epsilon} &=& \dfrac{1}{V_0A_j^2}\bigg[\pder{\nu}{\epsilon}\left(1+(\mu-\epsilon^2 \Btj)\dfrac{\tau-\tau_j}{2}\right)\\ &&+\tanh{\dfrac{\nu(\tau-\tau_j)}{2}}\left(\nu\pder{\nu}{\epsilon}\dfrac{\tau-\tau_j}{2}+\pder{\mu}{\epsilon}-2\epsilon\Btj-\epsilon^2\pder{\Btj}{\epsilon} \right)\bigg]. \\
\end{array}
$$
\subsection{Partial derivative of \texorpdfstring{$\psi$}{psi} with respect to \texorpdfstring{$\rho$}{rho}}
The partial derivative of $\psi$ with respect to $\rho$ writes: 
$$
    \dfrac{\partial \psi(X_{m,n}(t),V(t),t;z)}{\partial \rho} = \left[\dfrac{\partial \At}{\partial \rho} + V_0\dfrac{\partial \Bt}{\partial \rho} \right]\psi(X_{m,n}(t),V(t),t;z)
$$
with
$$
\def\arraystretch{2.7}
\begin{array}{rcl}
    \pder{\At}{\rho} &=& \pder{\Atj}{\rho}- \dfrac{\kappa\theta\lambda\tilde{\rho} z (\tau-\tau_j)}{\epsilon\rho}  + \dfrac{2\kappa\theta}{\epsilon^2}\pder{D_j(\tau)}{\rho},\\
    \pder{\Bt}{\rho} &=& \pder{\Btj}{\rho} - \dfrac{1}{V_0}\pder{A_j(\tau)}{\rho},
\end{array}
$$
where
$$
\def\arraystretch{2.7}
\begin{array}{rcl}
\pder{\mu}{\rho} &=& \dfrac{\mu-\kappa}{\rho}, \\
\pder{\nu}{\rho} &=& \dfrac{\mu}{\nu}\pder{\mu}{\rho},\\
\pder{D_j(\tau)}{\rho} &=& \dfrac{1}{\nu}\pder{\nu}{\rho} + \dfrac{1}{2}\left(\kappa\dfrac{\xi-1}{\rho}-\pder{\nu}{\rho}\right)(\tau-\tau_j)-\dfrac{1}{E_j(\tau)}\pder{E_j(\tau)}{\rho},  \\
\dfrac{1}{E_j(\tau)}\pder{E_j(\tau)}{\rho} &=& \dfrac{1}{V_0A_j^2(\tau)}\left[\pder{\nu}{\rho} +\left(\pder{\mu}{\rho}-\epsilon^2\pder{\Btj}{\rho}\right)\tanh{\dfrac{\nu(\tau-\tau_j)}{2}}\right] \\
&& - \dfrac{1}{E_j(\tau)}\pder{\nu}{\rho}(\tau-\tau_j)(\nu-\mu+\epsilon^2\Btj)e^{-\nu(\tau-\tau_j)}, \\
\pder{A_j(\tau)}{\rho} &=& \dfrac{1}{\tilde{A}_j^2(\tau)}\pder{\tilde{A}_j^1(\tau)}{\rho} - \dfrac{A_j(\tau)}{\tilde{A}_j^2(\tau)}\pder{\tilde{A}_j^2(\tau)}{\rho},\\
\dfrac{1}{\tilde{A}_j^2(\tau)}\pder{\tilde{A}_j^1(\tau)}{\rho} &=& \dfrac{1}{A_j^2(\tau)}\bigg[2\tanh{\dfrac{\nu(\tau-\tau_j)}{2}}\left(\mu\pder{\Btj}{\rho}-\epsilon^2\Btj\pder{\Btj}{\rho}+\Btj\pder{\mu}{\rho} \right)\\
&&+\pder{\nu}{\rho}\dfrac{\tau-\tau_j}{2}\left(\Btj(2\mu-\epsilon^2\Btj)+\lambda^2(z-z^2)\right)\bigg],\\
\dfrac{1}{\tilde{A}_j^2(\tau)}\pder{\tilde{A}_j^2(\tau)}{\rho} &=& \dfrac{1}{V_0A_j^2(\tau)}\bigg[\left(1+(\mu-\epsilon^2\Btj)\dfrac{\tau-\tau_j}{2}\right)\pder{\nu}{\rho} \\
&& + \left(\nu\pder{\nu}{\rho}\dfrac{t-T_j}{2} + \pder{\mu}{\rho} -\epsilon^2\pder{\Btj}{\rho}\right)\tanh{\dfrac{\nu(\tau-\tau_j)}{2}} \bigg].
\end{array}
$$
\subsection{Partial derivatives of \texorpdfstring{$\psi$}{psi} with respect to \texorpdfstring{$a$}{a}, \texorpdfstring{$b$}{b}, \texorpdfstring{$c$}{c} and \texorpdfstring{$d$}{d}}
One can observe that only $\gamma_k(\tau)$ depends on the parameters $a$, $b$, $c$ et $d$, which means that the derivatives of the characteristic function with respect to these parameters are close from each other. Consequently, we group the four partial derivatives in this section. Let $x \in \{a,b,c,d\}$, the partial derivative of $\psi$ with respect to $x$ writes 
$$
    \dfrac{\partial \psi(X_{m,n}(t),V(t),t;z) }{\partial x}= \left[\dfrac{\partial \At}{\partial x} + V_0\dfrac{\partial \Bt}{\partial x} \right]\psi(X_{m,n}(t),V(t),t;z)
$$
with
$$
\def\arraystretch{2.7}
\begin{array}{rcl}
    \pder{\At}{x} &=& \pder{\Atj}{x}- \dfrac{\kappa\theta z (\tau-\tau_j)}{\epsilon}\pder{(\tilde{\rho}\lambda)}{x}  + \dfrac{2\kappa\theta}{\epsilon^2}\pder{D_j(\tau)}{x},\\
    \pder{\Bt}{x} &=& \pder{\Btj}{x} - \dfrac{1}{V_0}\pder{A_j(\tau)}{x},
\end{array}
$$
where
$$
\def\arraystretch{2.7}
\begin{array}{rcl}
     \pder{(\tilde{\rho}\lambda)}{x} &=& \displaystyle\sum_{k=m}^{n-1}w_k(0)\pder{\|\gamma_k(\tau) \|}{x}\rho_k(\tau), \\
     \pder{\xi}{x} &=& \dfrac{\epsilon}{\kappa}\displaystyle\sum_{k=m}^{n-1}\alpha_k(0)\sum_{l=m(\tau)}^k\dfrac{\Delta T_l(F_l(0)+\delta)\rho_l(\tau)\pder{\|\gamma_l(\tau) \|}{x}}{1+\Delta T_l F_l(0)},\\
     \pder{\mu}{x} &=& \kappa \pder{\xi}{x} - \epsilon z\pder{(\tilde{\rho}\lambda)}{x} ,\\
     \pder{\lambda^2}{x} &=& 2\left\langle \displaystyle\sum_{k=m}^{n-1}w_k(0)\pder{\gamma_k(\tau) }{x},\displaystyle\sum_{k=m}^{n-1}w_k(0)\gamma_k(\tau) \right\rangle ,\\
     \pder{\nu}{x} &=& \dfrac{1}{\nu}\left(\pder{\mu}{x}\mu + \dfrac{1}{2}\pder{\lambda^2}{x}\epsilon^2(z-z^2) \right), \\
     \pder{D_j(\tau)}{x} &=& \dfrac{1}{\nu}\pder{\nu}{x} + \dfrac{1}{2}\left(\kappa\pder{\xi}{x} - \pder{\nu}{x}\right)(\tau-\tau_j) -\dfrac{1}{E_j(\tau)}\pder{E_j(\tau)}{x} ,\\
     \dfrac{1}{E_j(\tau)}\pder{E_j(\tau)}{x} &=& \dfrac{1}{V_0A_j^2(\tau)}\left[\pder{\nu}{x} +\left(\pder{\mu}{x}-\epsilon^2\pder{\Btj}{x}\right)\tanh{\dfrac{\nu(\tau-\tau_j)}{2}}\right] \\
     && - \dfrac{1}{E_j(\tau)}\pder{\nu}{x}(\tau-\tau_j)(\nu-\mu+\epsilon^2\Btj)e^{-\nu(\tau-\tau_j)},\\
     \pder{A_j(\tau)}{x} &=& \dfrac{1}{\tilde{A}_j^2(\tau)}\pder{\tilde{A}_j^1(\tau)}{x} - \dfrac{A_j(\tau)}{\tilde{A}_j^2(\tau)}\pder{\tilde{A}_j^2(\tau)}{x},\\
     \dfrac{1}{\tilde{A}_j^2(\tau)}\pder{\tilde{A}_j^1(\tau)}{x}&=& \dfrac{1}{A^2_j(\tau)}\bigg[2\tanh{\dfrac{\nu(\tau-\tau_j)}{2}} \bigg(\mu\pder{\Btj}{x}-\epsilon^2\Btj\pder{\Btj}{x}+\Btj\pder{\mu}{x}\\
     &&+\dfrac{1}{2}\pder{\lambda^2}{x}(z-z^2) \bigg)
     +\pder{\nu}{x}\dfrac{\tau-\tau_j}{2}\left(\Btj(2\mu-\epsilon^2\Btj)+\lambda^2(z-z^2)\right)\bigg],\\
     \dfrac{1}{\tilde{A}_j^2(\tau)}\pder{\tilde{A}_j^2(\tau)}{x} &=& \dfrac{1}{V_0A_j^2(\tau)}\bigg[\left(1+(\mu-\epsilon^2\Btj)\dfrac{\tau-\tau_j}{2}\right)\pder{\nu}{x} \\
     &&+ \left(\nu\pder{\nu}{x}\dfrac{\tau-\tau_j}{2} + \pder{\mu}{x} -\epsilon^2\pder{\Btj}{x}\right)\tanh{\dfrac{\nu(\tau-\tau_j)}{2}} \bigg].
     
\end{array}
$$

\section{On the Levenberg-Marquardt algorithm}
In this section, we detail the classic Levenberg-Marquardt algorithm and the extended version handling bound constraints. 
\subsection{Standard Levenberg-Marquardt algorithm} 
\label{classic_lm}
\begin{algorithm}[H] 
    \KwIn{$\BTheta_0$, $F$, $\f$, $\J$, L, $\epsilon_1$, $\epsilon_2$, $\epsilon_3$, $k_{max}$, $\tau$}  
    \Begin{
        $k \gets 0$; $\nu \gets 2$ \\
        $\mathbf{A} \gets \J(\BTheta_0)^T\J(\BTheta_0)$;  $\g \gets \J(\BTheta_0)^T\f(\BTheta_0)$ \\
        \textit{found} $\gets (F(\BTheta_k) \le \epsilon_1 \text{ or } \left\|\g\right\|_{\infty} \le \epsilon_2)$; $\mu \gets \tau  \max_i\{a_{ii}\}$ \\
        \While{!(\textit{found}) and $k < k_{max}$}{
            Solve $(\mathbf{A}+\mu\mathbf{I})\Bd = -\g$ \\
            \If{$\left\|\Bd\right\|^2 \le \epsilon_3^2\left\|\BTheta_k\right\|^2 $}{
                \textit{found} $\gets$ \textbf{true}
            }
            \Else{
                $\BTheta_{k+1} \gets \BTheta_k + \Bd$ \\
                \If{$F(\BTheta_k) - F(\BTheta_{k+1}) > 0$ and $L(\mathbf{0}) - L(\Bd) > 0$}{
                    $\eta \gets (F(\BTheta_k)-F(\BTheta_{k+1}))/(L(\mathbf{0})-L(\Bd))$\\
                    $\mathbf{A} \gets \J(\BTheta_{k+1})^T\J(\BTheta_{k+1})$; $\g \gets \J(\BTheta_{k+1})^T\f(\BTheta_{k+1})$ \\
                    \textit{trouve} $\gets (F(\BTheta_{k+1}) \le \epsilon_1 \text{ ou } \left\|\g\right\|_{\infty} \le \epsilon_2)$ \\
                    $\mu \gets \mu  \max\{\frac{1}{3},1-(2\eta-1)^3\}; \nu \gets 2$
                }
                \Else{
                    $\mu \gets \mu  \nu$; $\nu \gets 2\nu$ 
                }
            }
            $k \gets k +1 $
        }
    }
\caption{Levenberg-Marquardt algorithm}
\label{alg:lm_algo}
\end{algorithm}
In Algorithm \ref{alg:lm_algo}, the function $L$ corresponds to the value of the objective function $F$ in $\BTheta_{k+1}$ when the residuals are approximated by a first order Taylor expansion. Mathematically, we have:
\begin{equation*}
F(\BTheta_k + \Bd) \simeq L(\Bd) = \frac{1}{2}\left\|\f(\BTheta_k)+\J(\BTheta_k)\Bd \right\|^2
\end{equation*}
Hence, $L(\mathbf{0})-L(\Bd)$ can be interpreted as the gain predicted by a linear model. It is easy to check that this quantity is always positive. \\

The quantity $\eta$ (appearing first in line 13 of the routine above) allows to measure how good the approximation of $F(\BTheta_k + \Bd)$ by $L(\Bd)$ is. A large value of $\eta$ indicates that $L(\Bd)$ is a good approximation of $F(\BTheta_k + \Bd)$, whereas a small value of $\eta$ indicates the contrary. In the first case, $\mu$ is decreased in order to imitate the Gauss-Newton algorithm behaviour; in the second case, $\mu$ is increased in order to imitate the behaviour of the steepest descent method.

As regards the updating strategy of the damped parameter, we use the one introduced by \cite{NIELSEN1999}. 

\subsection{Extended Levenberg-Marquardt algorithm handling bound constraints} \label{extended_lm} 
We present an extension of the classic Levenberg-Marquardt algorithm which can handle bound constraints. This extension was first proposed by \cite{KANZOW2002}.

Let us denote by $P_X$ the projection onto the feasible set $X$ (which in the framework of the DDSVLMM is equal to $(\mathbb{R}_+)^4\times(\mathbb{R}_+^*)^3 \times ]-1;1[$). 
With respect to the Algorithm \ref{alg:lm_algo}, only the lines from 11 to 20 must be modified. They must be replaced by the following ones.

\begin{algorithm}[H]
        $\BTheta_{k+1} \gets P_X(\BTheta_k + \Bd)$; $\Bd \gets \BTheta_{k+1} - \BTheta_k$\\
        \If{$F(P_X(\BTheta_{k+1})) \le \gamma F(\Theta_k)$}{
            \If{$F(\BTheta_k) - F(\BTheta_{k+1}) > 0$ and $L(\mathbf{0}) - L(\Bd) > 0$}{
                $\eta \gets (F(\BTheta_k)-F(\BTheta_{k+1}))/(L(\mathbf{0})-L(\Bd))$\\
                $\mathbf{A} \gets \J(\BTheta_{k+1})^T\J(\BTheta_{k+1})$; $\g \gets \J(\BTheta_{k+1})^T\f(\BTheta_{k+1})$ \\
                \textit{found} $\gets (F(\BTheta_{k+1}) \le \epsilon_1 \text{ ou } \left\|\g\right\|_{\infty} \le \epsilon_2)$ \\
                $\mu \gets \mu  \max\{\frac{1}{3},1-(2\eta-1)^3\}; \nu \gets 2$
            }
            \Else{
                $\mu \gets \mu  \nu$; $\nu \gets 2\nu$ 
            }
        }
        \ElseIf{$\nabla F(\BTheta_{k+1})^T\Bd \le 0$}{
            Perform a line search, i.e. look for $\alpha$ such that $F(P_X(\BTheta_k+\alpha\Bd))$ is reasonably lower than $F(\BTheta_k)$
        }
        \Else{
            Apply a projected gradient step, i.e. compute $t = \max_{l\in \mathbb{N}} \beta^l$ such that $F(P_X(\BTheta_k-t\g)) \le F(\BTheta_k)+\sigma\nabla\g^T(P_X(\BTheta_k-t\g)-\BTheta_k)$ 
        }
        \caption{Extended Levenberg-Marquardt algorithm}
        \label{alg:lmbc_algo}
\end{algorithm}
The parameters $\gamma$, $\beta$ and $\sigma$ are empirically fixed parameters in $(0,1)$. 

\section{Results for USD market data}
\subsection{Methods accuracy}
\label{US_data_results}
\begin{figure}[H]
    \centering
    \includegraphics[width=\textwidth]{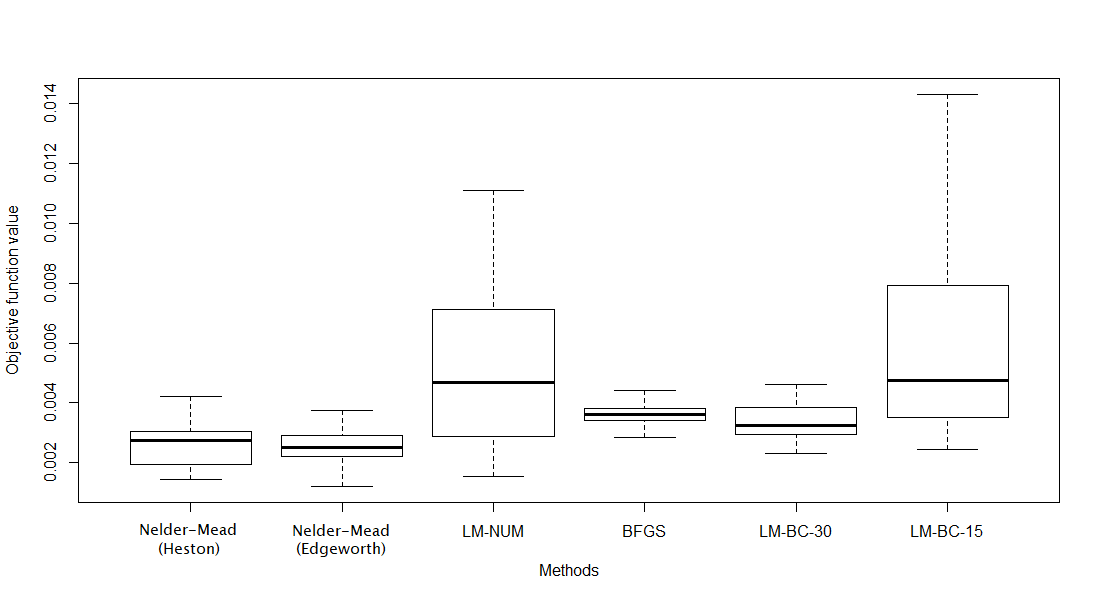}
    \caption{Boxplots comparing calibration methods of Table \ref{tab:summary_methods}}
    \label{fig:boxplots_usd}
\end{figure}

\begin{table}[H]
\centering
\begin{tabular}{|cc|}
\hline
Method &  \begin{tabular}[c]{@{}c@{}}Percentage of unsatisfied \\ Feller condition \end{tabular} \\ \hline
Heston & 0 \% \\
Edgeworth & 0 \% \\
LM-NUM & 25 \% \\
BFGS & 3 \%  \\
LM-BC-30 & 23 \% \\
LM-BC-15 & 23 \%  \\ \hline
\end{tabular}
\caption{Percentages of unsatisfied Feller condition over the 100 calibrations for each methods}
\label{tab:feller_condition_usd}
\end{table}

\begin{figure}[H]
    \centering
    \includegraphics{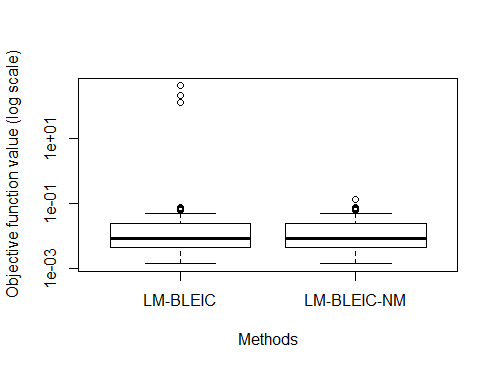}
    \caption{Boxplots comparing LM-BLEIC and LM-BLEIC-NM}
    \label{fig:boxplots_lmbleic_usd}
\end{figure}

\subsection{Time efficiency}
\label{US_time}
\begin{table}[H]
\resizebox{\textwidth}{!}{%
\begin{tabular}{|cccccc|}
\hline
Method & Average CPU time & \begin{tabular}[c]{@{}c@{}}Average number of \\ calls to $F$/$f$\end{tabular} & \begin{tabular}[c]{@{}c@{}}Average call time to \\ $F$/$f$\end{tabular} & \begin{tabular}[c]{@{}c@{}}Average number of \\ calls to $\nabla F$/$\nabla f$\end{tabular} & \begin{tabular}[c]{@{}c@{}}Average call time to\\ $\nabla F$/$\nabla f$\end{tabular} \\ \hline
Heston & 166.98 s & 1492.42 & 0.11 s & 0 & 0 s \\
Edgeworth & 9.00 s & 1487.06 & 5.75 $\times 10^{-3}$ s & 0 & 0 s \\
LM-NUM-FWD & 25.58 s & 177.55 & 0.14 s & 0 & 0 s \\
LM-NUM-CENT & 40.37 s & 282.85 & 0.14 s & 0 & 0 s \\
BFGS & 34.90 s & 37.89 & 0.14 s & 37.89 & 0.77 s \\
LM-BC-15 & 15.67 s & 42.70 & 0.14 s & 15.00 & 0.63 s \\
LM-BC-30 & 31.33 s & 85.06 & 0.14 s & 30.00 & 0.63 s \\
LM-BLEIC & 45.23 s & 98.61 & 0.14 s & 48.21 & 0.64 s \\
LM-BLEIC-NM & 47.21 s & 104.64 & 0.14 s & 48.21 & 0.64 s \\ \hline
\end{tabular}%
}
\caption{Computational times}
\label{tab:cpu_time_usd}
\end{table}

\end{document}